\documentclass[a4paper,10pt,fleqn]{amsart}
\usepackage{amssymb}
\usepackage{amsmath}
\usepackage{a4wide}
\usepackage{mathrsfs}

\newtheorem{thm}{Theorem}[section]
\newtheorem{lemma}[thm]{Lemma}
\newtheorem{ass}[thm]{Assumption}

\newtheorem{df}[thm]{Definition}

\theoremstyle{remark}
\newtheorem{rem}[thm]{Remark}
\newtheorem{exmp}[thm]{Example}
\newcommand{\MOG}{\mbox{Mogul'ski\u{\i}}}

\numberwithin{equation}{section}

\newcommand{\ee}[1]{\mathbb{E}\left[#1\right]}
\newcommand{\pp}[1]{\mathbb{P}\left(#1\right)}
\newcommand{\N}{\mathbb{N}}
\newcommand{\R}{\mathbb{R}}
\newcommand{\F}{\mathcal{F}}
\newcommand{\X}{\mathcal{X}}

\newcommand{\B}{\mathcal{B}}
\newcommand{\dd}{\mathrm{d}}
\newcommand{\beginproof}{\begin{proof}[\textbf{\emph{Proof}}]}
\newcommand{\beginproofof}[1]{\begin{proof}[\textbf{\emph{Proof of #1}}]}

\newcommand{\exx}[1]{\exp\left(#1\right)}
\newcommand{\lgg}[1]{\log\left(#1\right)}
\newcommand{\order}{\left(1+O\left(\frac{1}{n}\right)\right)}
\newcommand{\nt}{\notag}
\renewcommand{\phi}{\varphi}
\renewcommand{\epsilon}{\varepsilon}
\newcounter{todoCtr} \numberwithin{todoCtr}{section} 

\date{\today}

\author{V.\ Leijdekker}
\address{Korteweg-de Vries Institute for Mathematics,
University of Amsterdam, the Netherlands \: --- \: ABN A{\sc mro},
Amsterdam, the Netherlands} \curraddr{}
\email{v.j.g.leijdekker@uva.nl}

\author{M.\ Mandjes}
\address{Korteweg-de Vries Institute for Mathematics,
University of Amsterdam, the Netherlands \: --- \: E{\sc urandom},
Eindhoven, the Netherlands \: --- \: CWI, Amsterdam, the
Netherlands} \curraddr{} \email{m.r.h.mandjes@uva.nl}

\author{P.\ Spreij}
\address{Korteweg-de Vries Institute for Mathematics,
University of Amsterdam, the Netherlands} \curraddr{}
\email{p.j.c.spreij@uva.nl}

\title{Sample-path Large Deviations in Credit Risk}

\begin{document}

\begin{abstract}
The event of large losses plays an important role in credit risk.
As these large losses are typically rare, and portfolios usually
consist of a large number of positions, large deviation theory is
the natural tool to analyze the tail asymptotics of the
probabilities involved. We first derive a sample-path large
deviation principle (LDP) for the portfolio's loss process, which
enables the computation of the logarithmic decay rate of the
probabilities of interest. In addition, we derive exact asymptotic
results for a number of specific rare-event probabilities, such as
the probability of the loss process exceeding some given function.
\end{abstract}

\maketitle

\section{Introduction}
For financial institutions, such as banks and insurance companies,
it is of crucial  importance to accurately assess the risk of
their portfolios. These portfolios typically consist of a large
number of obligors, such as mortgages, loans or insurance
policies, and therefore it is computationally infeasible to treat
each individual  object in the portfolio separately. As a result,
attention has shifted to measures that characterize the risk of
the portfolio as a whole, see e.g.~\cite{book:DuffieSingleton} for
general principles concerning managing credit risk.  The best
known metric is the so-called {\it value at risk}, see
\cite{book:JorionVar}, measuring the minimum amount of money that
can be lost with $\alpha$ percent certainty over some given
period. Several other measures have been proposed, such as
economic capital, the risk adjusted return on capital ({\sc
raroc}), or expected shortfall, which is a coherent risk
measure~\cite{article:Artzneretal}. Each of these measure are
applicable to market risk as well as credit risk. Measures such as
loss given default ({\sc lgd}) and exposure at default ({\sc ead})
are measures that purely apply to credit risk. These and other
measures are discussed in detail in  e.g.
\cite{book:McNeilFreyEmbrechts}.

The currently existing methods  mainly focus on the distribution
of the portfolio loss {\it up to a  given point in time} (for
instance one year into the future). It can be argued, however,
that in many situations it makes more sense to use probabilities
that involve the (cumulative) loss {\it process}, say $\{L(t):
t\ge 0\}$. Highly relevant, for instance, is the event that
$L(\cdot)$ ever exceeds a given function $\zeta(\cdot)$ (within a
certain time window, for instance between now and one  year
ahead), i.e., an event of the type
\begin{equation}
\label{eq:event} \{\exists t \le T: L(t) \ge \zeta(t)\}.
\end{equation}
It is clear that measures of the latter type are intrinsically
harder to analyze, as it does not suffice anymore to have
knowledge of the marginal  distribution of the loss process at a
given point in time; for instance the event (\ref{eq:event})
actually corresponds to the union of events $\{L(t)\ge
\zeta(t)\}$, for $t\le T$, and its probability will depend on the
law of $L(\cdot)$ as a process on $[0,T]$.

In line with the remarks we made above, earlier papers on
applications of large-deviation theory to credit risk,  mainly
address the (asymptotics of the) distribution of the loss process
at a single point in time, see e.g.\
\cite{article:DemboDeuschelDuffie2004} and
\cite{article:GlassermanKangShahabuddin2007}. The former paper
considers, in addition, also the  probability that the {\it
increments} of the loss process exceed a certain level. Other
approaches to quantifying the tail distribution of the losses have
been taken by \cite{article:LucasKlaassenSpreijStraetmans2001},
who use extreme-value theory (see \cite{book:Embrechtsetal} for a
background), \cite{article:Gordy2002} and
\cite{article:MartinThompsonBrowne2001}, where the authors
consider saddle point approximations to the tails of the loss
distribution. Numerical and simulation techniques for credit risk
can be found in e.g.~\cite{book:Glasserman}. The first
contribution of our work concerns a so-called sample-path large
deviation principle (LDP) for the  average cumulative losses for
large portfolios. Loosely speaking, such an LDP means that, with
$L_n(\cdot)$ denoting the loss process when $n$ obligors are
involved, we can compute the  logarithmic asymptotics (for $n$
large) of the normalized loss process $L_n(\cdot)/n$ being in a
set of trajectories $A$:
\begin{equation}
\label{eq:logass} \lim_{n\to\infty}\frac{1}{n}\log {\mathbb
P}\left(\frac{1}{n}L_n(\cdot) \in A\right);
\end{equation}
we could for instance pick a set $A$ that  corresponds to the
event (\ref{eq:event}). Most of the sample-path LDPs that have
been developed so far  involve stochastic processes with
independent or nearly-independent increments, see for instance the
results by $\MOG$ for random walks \cite{article:Mogulskii1976},
de Acosta for L\'evy processes \cite{article:deacosta1994}, and
Chang \cite{article:chang1995} for weakly-correlated processes;
results for processes with a stronger correlation structure are
restricted to special classes of processes, such as Gaussian
processes, see e.g.\ \cite{article:Azencott1980}. It is observed
that our loss process is not  covered by these results, and
therefore new theory had to be developed. The proof of our LDP
relies on `classical' large-deviation results (such as Cram\'er's
theorem, Sanov's theorem, $\MOG$'s theorem), but in addition the
concept of {\it epi-convergence} \cite{article:PeterKall1986} is
relied upon.

Our second main result focuses specifically on the event
(\ref{eq:event}) of ever (before some time horizon $T$) exceeding
a given barrier function $\zeta(\cdot)$. Whereas we so far
considered,  inherently imprecise, logarithmic asymptotics of the
type displayed in (\ref{eq:logass}), we can now compute so-called
{\it exact asymptotics}: we identify an explicit function $f(n)$
such that $f(n)/p_n\to 1$ as $n\to \infty$, where $p_n$ is the
probability of our interest. As is known from the literature, it
is in general substantially harder to find exact asymptotics than
logarithmic asymptotics. The proof uses the fact that, after
discretizing time, the contribution of just a single time epoch
dominates, in the sense that there is a $t^\star$ such that
\begin{equation}
\label{eq:pn} \left.{\mathbb P}\left(\frac{1}{n} \,L_n(t^\star)\ge
\zeta(t^\star)\right)\right/p_n\to 1,\:\:\:\mbox{with}\:\:
p_n:={\mathbb P}\left(\exists t: \frac{1}{n}\, L_n(t)\ge
\zeta(t)\right).
\end{equation}
This $t^\star$ can be interpreted as the most likely epoch of
exceeding $\zeta(\cdot).$

Turning back to the setting of credit risk, both of the results we
present are derived in a setup where all obligors in the portfolio
are i.i.d., in the sense that they behave independently and
stochastically identically. A third contribution of our work
concerns a discussion on how to extend our results to cases where
the obligors are {\it dependent} (meaning that they, in the
terminology of \cite{article:DemboDeuschelDuffie2004}, react to
the same `macro-environmental' variable, conditional upon which
they are independent again). We also treat the case of
obligor-heterogeneity: we show how to extend the results to the
situation of multiple classes of obligors.

The paper is structured as follows.  In Section
\ref{sec:notAndDef} we introduce the loss process and we describe
the scaling under which we work. We also recapitulate  a couple of
relevant large-deviation results. Our first main result, the
sample-path LDP for the cumulative loss process, is stated and
proved in Section \ref{sec:samplePathLDP}. Special attention is
paid to, easily-checkable, sufficient conditions under which this
result holds. As argued above, the LDP is a generally applicable
result, as it yields an expression for the decay rate of any
probability  that depends on the entire sample path. Then, in
Section \ref{sec:exactResults}, we derive the exact asymptotic
behavior of the probability that, at some point in time, the loss
exceeds a certain threshold, i.e., the asymptotics of $p_n$, as
defined in (\ref{eq:pn}). After this we derive a similar result
for the {\it increments} of the loss process. Eventually, in
Section \ref{sec:discussion}, we discuss a number of possible
extensions to the results we have presented. Special attention is
given to allowing dependence between obligors, and to different
classes of obligors each having its own specific distributional
properties.

\section{Notation and Definitions}
\label{sec:notAndDef} The portfolios of banks and insurance
companies are typically very large: they may consist of several
thousands of assets. It is therefore computationally impossible to
estimate the risks for each element, or obligor, in a portfolio.
This explains why one attempts to assess the aggregated losses
resulting from defaults, e.g. bankruptcies, failure to repay loans
or insurance claims, for the portfolio as a whole. The risk in the
portfolio is then measured through this (aggregate) loss process.
In the following sections we introduce the loss process and the
portfolio constituents more formally.

\subsection{Loss Process}
\label{subsec:lossProcess} Let $(\Omega,\F,\mathbb{P})$ be the
probability space on which all random variables below are defined.
We assume that the portfolio consists of $n$ obligors and we
denote the default time of obligor $i$ by $\tau_i$. Further we
write $U_i$ for the loss incurred on a default of obligor $i$. We
then define the cumulative loss process $L_n$ as
\begin{align}
\label{eq:defLoss} L_n(t) &:= \sum_{i=1}^n U_i Z_i(t),
\end{align}
where $Z_i(t)=1_{\{\tau_i \leq t\}}$ is the default indicator of
obligor $i$. We assume that the loss amounts $U_i\geq 0$ are
i.i.d., and that the default times $\tau_i\geq 0$ are i.i.d.\ as
well. In addition we assume that the loss amounts and the default
times are mutually independent. In the remainder of this paper $U$
and $Z(t)$ denote generic random variables with the same
distribution as the $U_i$ and $Z_i(t)$, respectively.

Throughout this paper we assume that the defaults only occur on
the time grid ${\mathbb N}$; in Section~\ref{sec:discussion} we
discuss how to deal with the default epochs taking continuous
values. In some cases we explicitly consider a finite time grid,
say $\{1,2,\ldots,N\}$. The extension of the results we derive to
a more general grid $\{0<t_1<t_2<\ldots<t_N\}$ is completely
trivial. The distribution of the default times, for each $j$, is
denoted by
\begin{align}
\label{eq:defLossDist}
p_j&:= \pp{\tau = j},\\
\label{eq:defLossDistCum} F_j&:= \pp{\tau \leq j} = \sum_{i=1}^j
p_i.
\end{align}

Given the distribution of the loss amounts $U_i$ and the default
times $\tau_i$, our goal is to investigate the loss process. Many
of the techniques that have been developed so far, first fix a
time $T$ (typically one year), and then stochastic properties of
the cumulative loss at time $T$, i.e., $L_n(T)$, are studied.
Measures such as value at risk and economic capital are examples
of these `one-dimensional' characteristics. Many interesting
measures, however, involve properties of the entire {\it path} of
the loss process rather than those of just one time epoch,
examples being the probability that that $L_n(\cdot)$ exceeds some
barrier function $\zeta(\cdot)$ for some $t$ smaller than the
horizon $T$, or the probability that (during a certain period) the
loss always stays above a certain level. The event corresponding
to the former probability might require the bank to attract more
capital, or worse, it might lead to the bankruptcy of this bank.
The event corresponding to the latter event might also lead to the
bankruptcy of the bank, as a long period of stress may have
substantial negative implications. We conclude that having a
handle on these probabilities is therefore a useful instrument
when assessing the risks involved in the bank's portfolios.

As mentioned above, the number of obligors $n$ in a  portfolio is
typically very large, thus prohibiting analyses based on the
specific properties of the individual obligors. Instead, it is
more natural to study the asymptotical behavior of the loss
process as $n\to\infty$. One could rely on a central-limit-theorem
based approach, but in this paper we focus on rare events, by
using the theory of large deviations.

In the following subsection we provide some background of
large-deviation theory, and we define a number of quantities that
are used in the remainder of this paper.

\subsection{Large Deviation Principle}
\label{subsec:generalLDPtheory} In this section we give a short
introduction to the theory of large deviations. Here, in an
abstract setting, the limiting behavior of a family of probability
measures $\{\mu_n\}$ on the Borel sets $\B$ of a metric space
$(\X,d)$ is studied, as $n\to \infty$. This behavior is referred
to as the Large Deviation Principle (LDP), and it is characterized
in terms of a {\it rate function}. The LDP states lower and upper
exponential bounds for the value that the measures $\mu_n$ assign
to sets in a topological space $\X$. Below we state the definition
of the rate function that has been taken from
\cite{book:DemboZeitouni}.

\begin{df}
\label{df:rateFct} A rate function is a lower semicontinuous
mapping  $I:\X \to [0,\infty]$, for all $\alpha\in [0,\infty)$ the
level set $\Psi_I(\alpha):=\{x\vert\ I(x)\leq \alpha\}$ is a
closed subset of $\X$. A good rate function is a rate function for
which all the level sets are compact subsets of $\X$.
\end{df}

With the definition of the rate function in mind we state the
large deviation principle for the sequence of measure $\{\mu_n\}$.

\begin{df}
\label{df:LDP} We say that $\{\mu_n\}$ satisfies the large
deviation principle with a rate function $I(\cdot)$ if
\begin{itemize}
\item[(i)] (Upper bound) for any closed set $F\subseteq\X$
\begin{align}
\label{eq:defLDPupper} \limsup_{n\to\infty}
\frac{1}{n}\log{\mu_n(F)} \leq -\inf_{x\in F}I(x).
\end{align}
\item[(ii)] (Lower bound) for any open set $G\subseteq\X$
\begin{align}
\label{eq:defLDPlower} \liminf_{n\to\infty}
\frac{1}{n}\log{\mu_n(G)} \geq -\inf_{x\in G}I(x).
\end{align}
\end{itemize}
We say that a family of random variables $X=\{X_n\}$, with values
in $\X$, satisfies an LDP with rate function $I_X(\cdot)$ iff the
laws $\{\mu_n^X\}$ satisfy an LDP with rate function $I_X$, where
$\mu_n^X$ is the law of $X_n$.
\end{df}

The so-called Fenchel-Legendre transform plays an important role
in expressions for the rate function. Let for an arbitrary random
variable $X$, the logarithmic moment generating function,
sometimes referred to as cumulant generating function, be given by
\begin{align}
\label{eq:defCumFct} \Lambda_X(\theta) &:=
\log{M_X(\theta)}=\log{\ee{e^{\theta X}}}\leq \infty,
\end{align}
for $\theta \in \R$. The Fenchel-Legendre transform
$\Lambda_X^\star$ of $\Lambda_X$ is then defined by
\begin{align}
\label{eq:defFenLeg} \Lambda_X^\star(x)&:=
\sup_{\theta}\left(\theta x - \Lambda_X(\theta)\right).
\end{align}
We sometimes say that $\Lambda^\star_X$ is the Fenchel-Legendre
transform of $X$.

The LDP from Definition \ref{df:LDP} provides upper and lower
bounds for the log-asymptotic behavior of measures $\mu_n$. In
case of the loss process (\ref{eq:defLoss}), {\it fixed} at some
time $t$, we can easily establish an LDP by an application of
Cram\'er's theorem (Theorem \ref{thm:Cramer}). This theorem yields
that the rate function is given by $\Lambda^\star_{UZ(t)}(\cdot)$,
where $\Lambda_{UZ(t)}^\star(\cdot)$ is the Fenchel-Legendre
transform of the random variable $UZ(t)$.

The results we present in this paper involve either
$\Lambda_U^\star(\cdot)$ (Section \ref{sec:samplePathLDP}), which
corresponds to i.i.d.\ loss amounts $U_i$ only, or
$\Lambda_{UZ(t)}^\star(\cdot)$ (Section \ref{sec:exactResults}),
which corresponds to those loss amounts up to time $t$. In the
following section we derive an LDP for the whole path of the loss
process, which can be considered as an extension of Cram\'er's
theorem.

\section{A Sample-Path Large Deviation Result}
\label{sec:samplePathLDP} In the previous section we have
introduced the large deviation principle.  In this section we
derive a sample-path LDP for the cumulative loss
process~(\ref{eq:defLoss}). We consider the exponential decay of
the probability that the path of the loss process $L_n(\cdot)$ is
in some set $A$, as the size $n$ of the portfolio tends to
infinity.

\subsection{Assumptions}
\label{subsec:LDPassumptions} In order to state a sample-path LDP
we need to define the topology that we work on. To this end we
define the space ${\mathscr{S}}$ of all nonnegative and
nondecreasing functions on $T_N=\{1,2,\ldots,N\}$,
\begin{align}
\nt {\mathscr{S}} &:= \{f:T_N\to \R_0^+\vert \ 0\leq f_i\leq
f_{i+1}\ \mbox{for } i < N\}.
\end{align}
This set is identified with the space $\R_\leq^N:=\{x\in
\R^N\vert\ 0\leq x_i\leq x_{i+1}\ \mbox{for } i< N\}$. The
topology on this space is the one induced by the supremum norm
\begin{align}
\nt \vert\vert f \vert\vert_\infty &= \max_{i=1,\ldots,N}\vert f_i
\vert.
\end{align}

As we work on a finite-dimensional space, the choice of the norm
is not important, as any other norm on ${\mathscr{S}}$ would
result in the same topology. We use the supremum norm as this is
convenient in some of the proofs in this section.

We identify the space of all probability measures on $T_N$ with
the simplex $\Phi$:
\begin{align}
\label{eq:defPhiN} \Phi:=\left\{\phi\in \R^N\left\vert\
\sum_{i=1}^N \phi_i = 1,\ \phi_i\geq 0 \mbox{ for } i\leq
N\right.\right\}.
\end{align}
For a given $\phi\in\Phi$ we denote the cumulative distribution
function by $\psi$, i.e.,
\begin{align}
\label{eq:defPsi} \psi_i &= \sum_{j=1}^i \phi_j,\ \ \mbox{for }
i\leq N;
\end{align}
note that $\psi\in{\mathscr{S}}$ and $\psi_N=1$.

Furthermore, we consider the loss amounts $U_i$ as introduced in
Section \ref{subsec:lossProcess}, a $\phi\in\Phi$ with cdf $\psi$,
and a sequence of $\phi^n\in\Phi$, each with cdf $\psi^n$, such
that $\phi^n\to\phi$ as $n\to\infty$, meaning that
$\phi_i^n\to\phi_i$ for all $i\leq N$. We define two families of
measures $(\mu_n)$ and $(\nu_n)$,
\begin{align}
\label{eq:defMun}
\mu_n(A)&:=\pp{\left(\frac{1}{n}\sum_{j=1}^{[n \psi_i]}U_j\right)_{i=1}^N\in A},\\
\label{eq:defNun} \nu_n(A)&:=\pp{\left(\frac{1}{n}\sum_{j=1}^{[n
\psi_i^n]}U_j\right)_{i=1}^N\in A},
\end{align}
where $A\in\B:=\B(\mathbb{R}^N)$ and $[x]:= \sup\left\{k\in N\vert
k \leq x\right\}$. Below we state an assumption under which the
main result in this section holds. This assumption refers to the
definition of {\em exponential equivalence}, which can be found in
Definition \ref{df:expEquiv}.

\begin{ass}
\label{ass:ExpEquiv} Assume that $\phi^n\to \phi$ and moreover
that the measures $\mu_n$ and $\nu_n$ as defined in
$(\ref{eq:defMun})$ and $(\ref{eq:defNun})$, respectively, are
exponentially equivalent.
\end{ass}

>From Assumption \ref{ass:ExpEquiv} we learn that the differences
between the two measures $\mu_n$ and $\nu_n$ go to zero at a
`superexponential' rate. In the next section, in Lemma
\ref{lem:suffExpEquiv}, we provide a  sufficient condition, that
is easy to check, under which this assumption holds.

\subsection{Main Result}
\label{subsec:samplePathLDPmainResult} The assumptions and
definitions in the previous  sections allow us to  state the main
result of this section. We show that the average loss process
satisfies a large deviation principle as in (\ref{df:LDP}). This
principle allows us to approximate a large variety of
probabilities related to the average loss process, such as the
probability that the loss process stays above a certain time
dependent level or the probability that the loss process exceeds a
certain level before some given point in time.

\begin{thm}
\label{thm:resultFiniteGrid} With $\Phi$ as in
$(\ref{eq:defPhiN})$ and under Assumption $\ref{ass:ExpEquiv}$,
the average loss process, $ L_n(\cdot)/n$ satisfies an LDP with
rate function $I_{U,p}$. Here, for $x\in \R_\leq^N$, $I_{U,p}$ is
given by
\begin{align}
\label{eq:defRateFctFiniteGrid}
I_{U,p}(x)&:=\inf_{\phi\in\Phi}\sum_{i=1}^N \varphi_i
\left(\lgg{\frac{\varphi_i}{p_i}} + \Lambda_U^\star
\left(\frac{\Delta x_i}{\phi_i}\right)\right),
\end{align}
with $\Delta x_i := x_i - x_{i-1}$ and $x_0:=0$.
\end{thm}

Observing the rate function for this sample path LDP, we see that
the effects of the default times $\tau_i$ and the loss amounts
$U_i$ are nicely decomposed into the two terms in the rate
function, one involving the distribution of the default epoch
$\tau$ (the `Sanov term', cf.\ \cite[Thm.\
6.2.10]{book:DemboZeitouni}), the other one involving the incurred
loss size $U$ (the `Cram\'er term', cf.\ \cite[Thm.\
2.2.3]{book:DemboZeitouni}. Observe that we recover Cram\'er's
theorem by considering a time grid {\it consisting of a single
time point}, which means that Theorem \ref{thm:resultFiniteGrid}
extends Cram\'er's  result. We also remark that, informally
speaking, the optimizing $\varphi\in\Phi$ in
(\ref{eq:defRateFctFiniteGrid}) can be interpreted as the `most
likely' distribution of the loss epoch, given that the path of
$L_n(\cdot)/n$ is close to $x$.

As a sanity check we calculate the value of the rate function
$I_{U,p}(x)$ for the 'average path' of $L_n(\cdot)/n$, given by
$x_j^\star~=~\ee{U} F_j$ for $j\leq N$, where $F_j$ is the
cumulative distribution of the default times as given in
(\ref{eq:defLossDistCum}); this path should give a rate function
equal to 0. To see this we first remark that clearly
$I_{U,p}(x)\geq 0$ for all $x$, since both the Sanov term and the
Cram\'er term are non-negative. This yields the following chain of
inequalities:
\begin{align}
\nt
\lefteqn{0 \leq I_{U,p}(x^\star)=\inf_{\phi\in \Phi}\sum_{i=1}^N \phi_i\left(\lgg{\frac{\phi_i}{p_i}} + \Lambda_U^\star\left(\frac{\ee{U}p_i}{\phi_i}\right)\right)}\\
\nt
&\stackrel{\phi=p}{\leq\hfill} \sum_{i=1}^N p_i\left(\lgg{\frac{p_i}{p_i}} + \Lambda_U^\star\left(\frac{\ee{U}p_i}{p_i}\right)\right)\\
\nt &= \sum_{i=1}^N p_i
\Lambda_U^\star\left(\ee{U}\right)=\Lambda_U^\star(\ee{U})=0,
\end{align}
where we have used that for $\ee{U}<\infty$, it always holds that
$\Lambda_U^\star(\ee{U})=0$ cf. \cite[Lemma
2.2.5]{book:DemboZeitouni}. The inequalities above thus show that,
if the `average path' $x^\star$ lies in the set of interest, then
the corresponding decay rate is 0, meaning that the probability of
interest decays subexponentially.

In the proof of Theorem \ref{thm:resultFiniteGrid} we use the
following lemma, which is related to the concept of {\it
epi-convergence}, extensively discussed in
\cite{article:PeterKall1986}.

\begin{lemma}
\label{lem:switchLimsupSup} Let $f_n,f:D\to \R$, with
$D\subset\R^m$ compact. Assume that for all $x\in D$ and for all
$x_n\to x$ in $D$ we have
\begin{equation}
\label{eq:hypothesis:switchLimsupSup} \limsup_{n\to\infty}f_n(x_n)
\leq f(x).
\end{equation}
Then we have
\begin{equation}
\nt \limsup_{n\to\infty}\sup_{x\in D}f_n(x) \leq \sup_{x\in
D}f(x).
\end{equation}
\end{lemma}
\beginproof
Let $f_n^\star=\sup_{x\in D} f_n(x)$, $f^\star=\sup_{x\in D}
f(x)$. Consider a subsequence $f_{n_k}^\star \to
\limsup_{n\to\infty} f_n^\star$. Let $\epsilon >0$ and choose
$x_{n_k}$ such that $f_{n_k}^\star < f_{n_k}(x_{n_k}) + \epsilon$
for all $k$. By the compactness of $D$, there exists a limit point
$x$ such that along a subsequence $x_{n_{k_j}}\to x$. By the
hypothesis (\ref{eq:hypothesis:switchLimsupSup}) we then have \[
\limsup_{n\to\infty} f_n^\star\leq \lim  f_{n_{k_j}} (x_{n_{k_j}})
+\epsilon = f(x)+ \epsilon \leq f^\star+\epsilon.\] Let
$\epsilon\downarrow 0$ to obtain the result.
\end{proof}

\beginproofof{Theorem \ref{thm:resultFiniteGrid}}
We start by establishing an identity from which we show both
bounds. We need to calculate the probability
\begin{align}
\nt \pp{\frac{1}{n} L_n(\cdot) \in A} &=
\pp{\left(\frac{1}{n}L_{n}(1),\ldots,\frac{1}{n}L_{n}(N)\right)\in
A},
\end{align}
for certain $A\in\B$. For each point $j$ on the time grid $T_N$ we
record by  the `default counter' $K_{n,j}\in\{0,\ldots,n\}$ the
number of defaults at time $j$:
\begin{align}
\nt K_{n,j}&:=\#\{ i\in\{1,\ldots,n\}\:\vert\: \tau_i=j\}.
\end{align}
These counters allow us to rewrite the probability to
\begin{align}
\nt \pp{\frac{1}{n} L_n(\cdot) \in A} &=
\ee{\left.\pp{\left(\frac{1}{n}\sum_{j=1}^{K_{n,1}}U_{(j)},
\ldots,
\frac{1}{n}\sum_{j=1}^{K_{n,1}+\ldots+K_{n,N}}U_{(j)}\right)\in
A\:\right\vert \: K_n}}
\\
\label{eq:step2:resultFiniteGrid}
&=\sum_{k_1+\ldots+k_N=n}\pp{K_{n,i}=k_i\ \mbox{for } i\leq
N}\times \pp{\left(\frac{1}{n} \sum_{j=1}^{m_i}
U_{(j)}\right)_{i=1}^N\in A},
\end{align}
where $m_i := \sum_{j=1}^i k_j$ and the loss amounts $U_j$ have
been ordered, such that the first $U_{(j)}$ correspond to the
losses at time 1, etc.

{\it Upper bound.} Starting from Equality
(\ref{eq:step2:resultFiniteGrid}), let us first establish the
upper bound of the LDP. \ To this end, let $F$ be a closed set and
consider the decay rate
\begin{align}
\label{eq:limsup1}
\limsup_{n\to\infty}\frac{1}{n}\log{\pp{\frac{1}{n} L_n(\cdot) \in
F}}.
\end{align}
An application of Lemma \ref{lem:app:DZ1.2.15} together with
(\ref{eq:step2:resultFiniteGrid}), implies that (\ref{eq:limsup1})
equals
\begin{align}
\nt
\lefteqn{\limsup_{n\to\infty}\frac{1}{n}\log{\pp{\frac{1}{n} L_n(\cdot) \in F}}}&\\
\label{eq:logSplit} &=\:
\limsup_{n\to\infty}\max_{\sum{k_j}={n}}\frac{1}{n}
\left[\log{\pp{\frac{K_{n,i}}{n}=\frac{k_i}{n}, \ i\leq
N}}+\log{\pp{\left(\frac{1}{n}\sum_{j=1}^{m_i}U_{(j)}\right)_{i=1}^N\in
F}}\right],
\end{align}
Next we replace the dependence on $n$ in the maximization  by
maximizing over the set $\Phi$ as in  (\ref{eq:defPhiN}). In
addition, we replace the $k_i$ in (\ref{eq:logSplit}) by
\begin{align}
\label{eq:defHatPhi} \hat\phi_{n,i}&:=
\frac{[n\psi_i]-[n\psi_{i-1}]}{n},
\end{align}
where the $\psi_i$ have been defined in $(\ref{eq:defPsi})$. As a
result, (\ref{eq:limsup1}) reads
\begin{align}
\label{eq:preMaxSwitch} &\:\limsup_{n\to\infty}\sup_{\varphi\in
\Phi}\frac{1}{n} \left[\log{\pp{\frac{K_{n,i}}{n}= \hat
\phi_{n,i}, \ i\leq N}}
 +\log{\pp{\left(\frac{1}{n}\sum_{j=1}^{[n\psi_i]}U_{(j)}\right)_{i=1}^N\in F}}\right].
\end{align}
Note that (\ref{eq:limsup1}) equals (\ref{eq:preMaxSwitch}), since
for each $n$ and vector $\left(k_1,\ldots,k_N\right)\in\N^N$, with
$\sum_{i=1}^Nk_i=n$, there is a $\phi\in \Phi$ with
$\phi_i=k_i/n$. On the other hand, we only cover outcomes of this
form by rounding off the $\phi_i$.

We can bound the first term in this expression from above using
Lemma \ref{lem:app:DZ2.1.9}, which implies that the decay rate
(\ref{eq:limsup1}) is majorized by
\begin{align}
\nt &\limsup_{n\to\infty}\sup_{\varphi\in \Phi}
\left[-\sum_{i=1}^N \hat
\phi_{n,i}\lgg{\frac{\hat\phi_{n,i}}{p_i}}
+\frac{1}{n}\log{\pp{\left(\frac{1}{n}
\sum_{j=1}^{[n\psi_i]}U_{(j)}\right)_{i=1}^N \in F}}\right].
\end{align}
Now note that calculating the limsup in the previous expression is
not straightforward due to the supremum over $\Phi$. The idea is
therefore to interchange the supremum and the limsup, by using
Lemma \ref{lem:switchLimsupSup}. To apply this lemma we first
introduce
\begin{align}
\nt
f_n(\phi)&:= -\sum_{i=1}^N \hat \phi_{n,i}\lgg{\frac{\hat\phi_{n,i}}{p_i}} +\frac{1}{n}\log{\pp{\left(\frac{1}{n}\sum_{j=1}^{[n\psi_i]}U_{(j)}\right)_{i=1}^N\in F}},\\
\nt f(\phi)&:=   -\inf_{x\in F}\sum_{i=1}^N
\phi_{i}\left(\lgg{\frac{\phi_{i}}{p_i}} + \Lambda_U^\star
\left(\frac{\Delta x_i}{\phi_i}\right)\right),
\end{align}
and note that $\Phi$ is a compact subset of $\R^n$. We have to
show that for any sequence $\phi^n\to \phi$
\begin{align}
\label{eq:targetInequalitySwitchLimsupSup}
\limsup_{n\to\infty}f_n(\phi^n)&\leq f(\phi),
\end{align}
such that the conditions of Lemma \ref{lem:switchLimsupSup} are
satisfied. We observe, with $\psi_i^n$ as in (\ref{eq:defPsi}) and
$\hat \phi_i^n$ as in (\ref{eq:defHatPhi}) with $\phi$ replaced by
$\phi^n$, that
\begin{align}
\nt
\limsup_{n\to\infty}f_n(\phi^n) &\leq \limsup_{n\to\infty}\left( -\sum_{i=1}^N \hat \phi_{n,i}^n\lgg{\frac{\hat\phi_{n,i}^n}{p_i}}\right)\\
\nt &\quad + \limsup_{n\to\infty}\frac{1}{n}
\log{\pp{\left(\frac{1}{n}
\sum_{j=1}^{[n\psi_i^n]}U_{(j)}\right)_{i=1}^N\in F}}.
\end{align}
Since $\phi^n\to \phi$ and since  $\hat \phi_{n,i}^n$ differs at
most by ${1}/{n}$ from $\phi_i^n$, it immediately follows that
$\hat\phi_{n,i}^n\to \phi_i$. For an arbitrary continuous function
$g$ we thus have $g\left(\hat\phi_{n,i}^n\right)\to g(\phi_i)$.
This implies that
\begin{align}
\label{eq:sanovUpper} \limsup_{n\to\infty}\left( -\sum_{i=1}^N
\hat \phi_{n,i}\lgg{\frac{\hat\phi_{n,i}^n}{p_i}}\right) &=
-\sum_{i=1}^N \phi_{i}\lgg{\frac{\phi_{i}}{p_i}}.
\end{align}
Inequality (\ref{eq:targetInequalitySwitchLimsupSup}) is
established once we have shown that
\begin{align}
\label{eq:MogUpper} \limsup_{n\to\infty}\frac{1}{n}
\log{\pp{\left(\frac{1}{n}
\sum_{j=1}^{[n\psi_i^n]}U_{(j)}\right)_{i=1}^N\in F}}\leq
-\inf_{x\in F}\sum_{i=1}^N\phi_i \Lambda_U^\star
\left(\frac{\Delta x_i}{\phi_i}\right).
\end{align}
By Assumption \ref{ass:ExpEquiv}, we can exploit the exponential
equivalence together with Theorem \ref{thm:app:DZ4.2.13}, to see
that (\ref{eq:MogUpper}) holds as soon as we have that
\begin{align}
\nt \limsup_{n\to\infty}\frac{1}{n} \log{\pp{\left(\frac{1}{n}
\sum_{j=1}^{[n\psi_i]}U_{(j)}\right)_{i=1}^N\in F}} \leq
-\inf_{x\in F}\sum_{i=1}^N\phi_i \Lambda_U^\star\left(\frac{\Delta
x_i}{\phi_i}\right).
\end{align}
But this inequality is a direct consequence of Lemma
\ref{lem:app:DZ5.1.8}, and we conclude that (\ref{eq:MogUpper})
holds. Combining (\ref{eq:sanovUpper}) with (\ref{eq:MogUpper})
yields
\begin{align}
\nt
\limsup_{n\to\infty}f_n(\phi^n) &\leq -\sum_{i=1}^N \phi_{i}\lgg{\frac{\phi_{i}}{p_i}} -\inf_{x\in F}\sum_{i=1}^N\phi_i \Lambda_U^\star\left(\frac{\Delta x_i}{\phi_i}\right)\\
\nt &= - \inf_{x\in F}\sum_{i=1}^N\phi_i
\left(\lgg{\frac{\phi_{i}}{p_i}} +
\Lambda_U^\star\left(\frac{\Delta
x_i}{\phi_i}\right)\right)=f(\phi),
\end{align}
so that indeed the conditions of Lemma \ref{lem:switchLimsupSup}
are satisfied, and therefore
\begin{align}
\nt
\limsup_{n\to\infty}\sup_{\phi\in \Phi} f_n(\phi)&\leq \sup_{\phi\in\Phi} f(\phi)= \sup_{\phi\in\Phi}\left(- \inf_{x\in F}\sum_{i=1}^N\phi_i \left(\lgg{\frac{\phi_{i}}{p_i}} + \Lambda_U^\star\left(\frac{\Delta x_i}{\phi_i}\right)\right)\right)\\
\nt &= -\inf_{x\in F}\inf_{\phi\in\Phi}\sum_{i=1}^N\phi_i
\left(\lgg{\frac{\phi_{i}}{p_i}} +
\Lambda_U^\star\left(\frac{\Delta x_i}{\phi_i}\right)\right)=
-\inf_{x\in F}I_{U,p}(x).
\end{align}
This establishes the upper bound of the LDP.

{\it Lower bound.} To complete the proof, we need to establish the
corresponding lower bound. Let $G$ be an open set and consider
\begin{align}
\label{liminfP} \liminf_{n\to\infty}
\frac{1}{n}\log{\pp{\frac{1}{n}L_n(\cdot)\in G}}.
\end{align}
We apply Equality (\ref{eq:step2:resultFiniteGrid}) to this
liminf, with $A$ replaced by $G$, and we observe that this sum is
larger than the largest term in the sum, which shows that (where
we directly switch to the enlarged space $\Phi$) the decay rate
(\ref{liminfP}) majorizes
\begin{align}
\nt &\liminf_{n\to\infty}\sup_{\varphi\in \Phi}\frac{1}{n}
\left(\log{\pp{\frac{1}{n}K_{n,i}= \hat\phi_{n,i}, \ i\leq N}}
+\log{\pp{\left(\frac{1}{n}\sum_{j=1}^{[n\psi_i]}U_{(j)}\right)_{i=1}^N\in
G}}\right).
\end{align}
Observe that for any sequence of functions $h_n(\cdot)$ it holds
that $\liminf_n \sup_x h_n(x) \geq \liminf_n h_n(\tilde x)$ for
all $\tilde x,$ so that we obtain the evident inequality
\begin{align}
\nt \liminf_{n\to\infty} \sup_x h_n(x) &\geq
\sup_{x}\liminf_{n\to\infty} h_n(x).
\end{align}
This observation yields that the decay rate of interest
(\ref{liminfP}) is not smaller than
\begin{align}
\label{eq:twoLiminfs} &\sup_{\varphi\in
\Phi}\left(\liminf_{n\to\infty}\frac{1}{n}
 \log{\pp{\frac{1}{n}K_{n,i}= \hat\phi_{n,i}, \ i\leq N}}+
 \liminf_{n\to\infty}\frac{1}{n}\log{\pp{\hspace{-1mm}\left(\frac{1}{n}
 \sum_{j=1}^{[n\psi_i]}U_{(j)}\right)_{i=1}^N\hspace{-2mm}\in G}}\hspace{-1mm}\right),
\end{align}
where we have used that $\liminf_n(x_n+y_n)\geq \liminf_n
x_n+\liminf_n y_n$. We apply Lemma \ref{lem:app:DZ2.1.9} to the
first liminf in (\ref{eq:twoLiminfs}), leading to
\begin{align}
\nt \lefteqn{
\liminf_{n\to\infty} \frac{1}{n} \log \pp{\frac{1}{n}K_{n,i}= \hat\phi_{n,i}, \ i\leq N}\geq \liminf_{n\to \infty}\left( -\sum_{i=1}^N \hat\phi_{n,i} \lgg{\frac{\hat\phi_{n,i}}{p_i}} -\frac{N}{n}\log(n+1)\right)}\\
\:&=\:- \sum_{i=1}^N \phi_i \log\left(\frac{\phi_i}{p_i}\right),
\end{align}
since $\log(n+1)/n\to 0$ as $n\to\infty.$ The second liminf in
(\ref{eq:twoLiminfs}) can be bounded from below by an application
of Lemma \ref{lem:app:DZ5.1.8}. Since $G$ is an open set, this
lemma yields
\begin{align}
\nt \liminf_{n\to\infty}\frac{1}{n}\log{\pp{\left(\frac{1}{n}
\sum_{j=1}^{[n\psi_i]}U_{(j)}\right)_{i=1}^N\in G}} \geq -
\inf_{x\in G}\sum_{i=1}^N\phi_i
\Lambda_U^\star\left(\frac{x_i-x_{i-1}}{\phi_i}\right).
\end{align}
Upon combining these two results, we see that we have established
the lower bound:
\begin{align}
\nt \lefteqn{\liminf_{n\to \infty} \frac{1}{n}\log{\pp{\frac{1}{n} L_n(\cdot) \in G}}}\\
\nt& \geq -\inf_{\phi\in\Phi} \inf_{x\in A}\left(\sum_{i=1}^N
\varphi_i\left(\log\left(\frac{\phi_i}{pi}\right) +
\Lambda_U^\star\left(\frac{x_i-x_{i-1}}{\phi_i}\right)\right)\right)=-\inf_{x\in
G}I_{U,p}(x).
\end{align}
This completes the proof of the theorem.
\end{proof}

In order to apply Theorem \ref{thm:resultFiniteGrid}, one needs to
check that Assumption \ref{ass:ExpEquiv} holds. In general this
could be a quite cumbersome exercise. In Lemma
\ref{lem:suffExpEquiv} below we provide a sufficient,
easy-to-check condition under which this assumption holds.

\begin{lemma}
\label{lem:suffExpEquiv} Assume that for all $\theta\in \R:$
$\Lambda_U(\theta)<\infty$. Then Assumption $\ref{ass:ExpEquiv}$
holds.
\end{lemma}

\begin{rem}
The assumption we make in Lemma \ref{lem:suffExpEquiv}, i.e., that
the logarithmic moment generating function is finite everywhere,
is a common assumption in large deviations theory. We remark that
for instance ${\MOG}$'s theorem \cite[Thm.\
5.1.2]{book:DemboZeitouni} also relies on this assumption; this
theorem is a sample-path LDP for
\[Y_n(t):=\frac{1}{n}\sum_{i=1}^{[nt]}X_i,\]
on the interval $[0,1]$. In $\MOG$'s result the $X_i$ are assumed
to be i.i.d; in our model we have that $L_n(t)=\sum_{i=1}^n
U_iZ_i(t)/n$, so that our sample-path result clearly does not fit
into the setup of $\MOG$'s theorem.\hfill$\diamondsuit$
\end{rem}

\begin{rem}
\label{rem:equivalence} In Lemma \ref{lem:suffExpEquiv} it was
assumed that $\Lambda_U(\theta)<\infty$, for all $\theta\in \R$,
but an equivalent condition is
\begin{align}
\label{eq:suffLambdaStar}
\lim_{x\to\infty}\frac{\Lambda_U^\star(x)}{x}=\infty.
\end{align}
In other words: this alternative condition can be used instead of
the condition stated in Lemma  \ref{lem:suffExpEquiv}. To see that
both requirements are equivalent, make the following observations.
In Lemma
\ref{lem:app:DZ2.2.20} it is shown that (\ref{eq:suffLambdaStar})
is implied by the assumption in Lemma \ref{lem:suffExpEquiv}. In
order to prove the converse, assume that (\ref{eq:suffLambdaStar})
holds, and that there is a $0<\theta_0<\infty$ for which
$\Lambda_U(\theta)=\infty$. Without loss of generality we can
assume that $\Lambda_U(\theta)$ is finite for $\theta <\theta_0$
and infinite for $\theta\geq\theta_0$. For $x>\ee{U}$, the
Fenchel-Legendre transform is then given by
\begin{align}
\nt \Lambda_U^\star(x)&= \sup_{0<\theta<\theta_0}\left(\theta x -
\Lambda_U(\theta)\right).
\end{align}
Since $U\geq 0$ and $\Lambda_U(0)=0$,  we know that
$\Lambda_U(\theta)\geq 0$ for $0<\theta<\theta_0$, and hence
\begin{align}
\nt \frac{\Lambda_U^\star(x)}{x}&\leq \theta_0,
\end{align}
which contradicts the assumption that this ratio tends to infinity
as $x\to\infty$, and thus establishing the
equivalence.\hfill$\diamondsuit$
\end{rem}

\beginproofof{Lemma \ref{lem:suffExpEquiv}}

Let $\phi^n\to \phi$ for some sequence of $\phi^n\in\Phi$ and
$\phi\in\Phi$. We introduce two families of random variables
$\{Y_n\}$ and $\{Z_n\}$,
\begin{align}
\nt Y_n &:= \left(\frac{1}{n}\sum_{j=1}^{[n
\psi_i]}U_j\right)_{i=1}^N,\:\:\:\: Z_n :=
\left(\frac{1}{n}\sum_{j=1}^{[n \psi_i^n]}U_j\right)_{i=1}^N,
\end{align}
which have laws $\mu_n$ and $\nu_n$, respectively, as in
(\ref{eq:defMun})--(\ref{eq:defNun}). Since $\phi^n\to \phi$ we
know that for any $\epsilon>0$ there exists an $M_\epsilon$ such
that for all $n>M_\epsilon$ we have that
$\max_i~\vert~\phi_i^n~-~\phi_i\vert~<~\epsilon/N$, and thus
$\vert \psi_i^n - \psi_i\vert < \epsilon$.

We have to show that for any $\delta>0$
\begin{align}
\nt \limsup_{n\to\infty} \frac{1}{n}\log{\pp{\vert \vert
Y_n-Z_n\vert \vert_\infty>\delta}} &= -\infty.
\end{align}
For $i\leq N$ consider the absolute difference between $Y_{n,i}$
and $Z_{n,i}$, i.e.,
\begin{align}
\label{eq:lem:expEquiv:step1} \left\vert Y_{n,i} - Z_{n,i}
\right\vert &= \left\vert \frac{1}{n}\sum_{j=1}^{[n\psi_i]}U_j -
\frac{1}{n}\sum_{j=1}^{[n\psi_i^n]}U_j\right\vert.
\end{align}
Next we have that for any $n>M_\epsilon$ it holds that $\vert
n\psi_i^n - n\psi_i\vert < n\epsilon,$ which yields the upper
bound,
\begin{align}
\nt \vert [n\psi_i^n]-[n\psi_i]\vert &< [n\epsilon] + 2,
\end{align}
since the rounded numbers differ at most by 1 from their real
counterparts. This means that the difference of the two sums in
(\ref{eq:lem:expEquiv:step1}) can be bounded by at most
$[n\epsilon] + 2$ elements of the $U_j$, which are for convenience
denoted by $U_j^\star$. Recalling that the $U_j$ are nonnegative,
we obtain
\begin{align}
\nt \left\vert \frac{1}{n}\sum_{j=1}^{[n \psi_i]}U_j -
\frac{1}{n}\sum_{j=1}^{[n \psi_i^n]}U_j\right\vert &\leq
\frac{1}{n}\sum_{j=1}^{[n\epsilon] + 2}U_j^\star.
\end{align}
Next we bound the probability that the difference exceeds
$\delta$, by using the above inequality:
\begin{align}
\nt \pp{\vert \vert Y_n-Z_n\vert \vert>\delta} \:&\leq\:
\pp{\frac{1}{n}\sum_{j=1}^{[n\epsilon] +
2}U_j^\star>\delta}\:\leq\: \ee{\exx{\theta U_1}}^{[n\epsilon] +
2}e^{-n\delta\theta},
\end{align}
where the last inequality follows from the Chernoff bound
\cite[Eqn.\ (2.2.12)]{book:DemboZeitouni} for $\theta>0$. Taking
the log of this probability, dividing by $n$, and taking the
limsup on both sides results in
\begin{align}
\nt \limsup_{n\to\infty} \frac{1}{n}\log{\pp{\vert \vert
Y_n-Z_n\vert \vert_\infty>\delta}} &\leq
\epsilon\Lambda_U(\theta)-\delta\theta.
\end{align}
By the assumption, $\Lambda_U(\theta)<\infty$ for all $\theta$.
Thus, $\epsilon\to 0$ yields
\begin{align}
\nt \limsup_{n\to\infty} \frac{1}{n}\log{\pp{\vert \vert
Y_n-Z_n\vert \vert_\infty>\delta}} &\leq -\delta\theta.
\end{align}
As $\theta$ was arbitrary, the exponential equivalence follows by
letting $\theta\to\infty$.
\end{proof}

We conclude this section with some examples.

\begin{exmp}
\label{exmp:finiteSupport} Assume that the loss amounts have
finite support, say on the interval $[0,u]$. Then we clearly have
\begin{align}
\nt \Lambda_U(\theta)&=\log{\ee{e^{\theta U}}} \leq \theta u <
\infty.
\end{align}
So for any distribution with finite support,  the assumption for
Lemma \ref{lem:suffExpEquiv} is satisfied, and thus Theorem
\ref{thm:resultFiniteGrid} holds. Here, the i.i.d.\ default times,
$\tau_i$, can have an arbitrary discrete distribution on the time
grid $\{1,\ldots,N\}$.

In practical applications  one (always) chooses a distribution
with finite support for the loss amounts, since the exposure to
every obligor is finite. Theorem \ref{thm:resultFiniteGrid} thus
clearly holds for any (realistic) model of the loss given default.

An explicit  expression for the rate function
(\ref{eq:defRateFctFiniteGrid}), or even the Fenchel-Legendre
transform, is usually not available. On the other hand one can use
numerical optimization techniques to calculate these
quantities.\hfill$\diamondsuit$
\end{exmp}

We next present an example to which Lemma~\ref{lem:suffExpEquiv}
applies.
\begin{exmp}
Assume that the loss amount $U$ is measured in a certain unit, and
takes on the values $u,2u,\ldots$ for some $u>0$. Assume that it
has a distribution of Poisson type with parameter $\lambda>0$, in
the sense that for $i=1,2,\ldots$
\begin{align}
\nt \pp{U=(i+1)u}&=e^{-\lambda}\frac{\lambda^i}{i!}.
\end{align}
It is then easy to check that $\Lambda_U(\theta)= \theta u +
\lambda\left(e^{\theta u}-1\right)$, being finite for all
$\theta$. Further calculations yield
\begin{align}
\nt \Lambda_U^\star(x)=
\left(\frac{x}{u}-1\right)\lgg{\frac{1}{\lambda}\left(\frac{x}{u}-1\right)}
- \left(\frac{x}{u}-1\right)+\lambda
\end{align}
for all $x>u$, and $\infty$ otherwise. Dividing this expression by
$x$ and letting $x\to\infty$, we observe that the resulting ratio
goes to $\infty$. As a consequence, Remark~\ref{rem:equivalence}
now entails that Theorem \ref{thm:resultFiniteGrid} applies. It
can also be argued that for any distribution $U$ with tail
behavior comparable to that of a Poisson distribution,
Theorem~\ref{thm:resultFiniteGrid} applies as well.
\hfill$\diamondsuit$
\end{exmp}

\section{Exact Asymptotic Results}
\label{sec:exactResults} In the previous section we have
established a sample-path large deviation principle on a finite
time grid; this LDP provides us with logarithmic asymptotics of
the probability that the sample path of $L_n(\cdot)/n$ is
contained in a given set, say $A$. The results presented in this
section are different in several ways. In the first place, we
derive exact asymptotics (rather than logarithmic asymptotics). In
the second place, our time domain is not assumed to be finite;
instead we consider all integer numbers, $\N$. The price to be
paid is that we restrict ourselves to special sets $A$, viz.\
those corresponding to the loss process (or the increment of the
loss process) exceeding a given function. We work under the setup
that we introduced in Section \ref{subsec:lossProcess}.

\subsection{Crossing a Barrier}
In this section we consider the asymptotic behavior of the
probability that the loss process at some point in time is above a
time-dependent level $\zeta$. More precisely, we consider the set
\begin{align}
\label{eq:crossingSet} A:= \left\{f:T\to\R_0^+\vert\ \exists t\in
T:\ f(t) \geq \zeta(t)\right\},
\end{align}
for some function $\zeta(t)$ satisfying
\begin{align}
\label{eq:interestingLevel} \zeta(t) >
\ee{UZ(t)}=\ee{U}F_t\:\:\:\:\mbox{for all $t\in T$,}
\end{align}
with $F_t$ as in (\ref{eq:defLossDistCum}). If we would consider a
function $\zeta$ that does not satisfy
(\ref{eq:interestingLevel}), we are not in a large-deviations
setting, in the sense that the probability of the event
$\{L_n(\cdot)/n \in A\}$ converges to 1 by the law of large
numbers. In order to obtain a more interesting result, we thus
limit ourselves to levels that satisfy
(\ref{eq:interestingLevel}). For such levels we state the first
main result of this section.

\begin{thm}
\label{thm:crossing} Assume that
\begin{align}
\label{eq:thm:crossing:unicity} \mbox{there is a unique
$t^\star\in T$ such that}\:\: I_{UZ}(t^\star)&=\min_{t\in T}
I_{UZ}(t),
\end{align}
and that
\begin{align}
\label{eq:thm:crossing:limit} \liminf_{t\to\infty}
\frac{I_{UZ}(t)}{\log t}&>0,
\end{align}
where $I_{UZ}(t) = \sup_{\theta}\left\{\theta
\zeta(t)-\Lambda_{UZ(t)}(\theta)\right\}=\Lambda_{UZ(t)}^\star\left(\zeta(t)\right).$
Then
\begin{align}
\label{eq:thm:crossing:result} \pp{\frac{1}{n}L_n(\cdot)\in A} &=
\frac{e^{-n I_{UZ}(t^\star)}C^\star}{\sqrt{ n }} \left( 1+O \left(
\frac{1}{n} \right)\right),
\end{align}
for $A$ as in $(\ref{eq:crossingSet})$ and $\sigma^\star$ is such
that $\Lambda_{UZ(t^\star)}'(\sigma^\star) =\zeta(t^\star)$. The
constant $C^\star$ follows from the Bahadur-Rao theorem (Theorem
\ref{thm:bahadurRao}), with
$C^\star=C_{UZ(t^\star),\zeta(t^\star)}$.
\end{thm}

Before proving our result, which will rely on arguments similar to
those in~\cite{article:LikhanovMazumdar1999}, we first discuss the
meaning and implications of Theorem \ref{thm:crossing}. In
addition we reflect on the role played by the assumptions. We do
so by a sequence of remarks.

\begin{rem}
Comparing Theorem \ref{thm:crossing} to the Bahadur-Rao theorem
(Theorem \ref{thm:bahadurRao}), we observe that the probability of
a sample mean exceeding a rare value has the same type of decay as
the probability of our interest (i.e., the probability that the
normalized loss process $L_n(\cdot)/n$ ever exceeds some function
$\zeta$). This decay looks like $C e^{-nI}/\sqrt{n}$ for positive
constants $C$ and $I$. This similarity can be explained as
follows.

First observe that the probability of our interest is actually the
probability of a {\em union} events. Evidently, this probability
is larger than the probability of any of the events in this union,
and hence also larger than the largest among these:
\begin{align}
\label{eq:unionBound} \pp{\frac{1}{n} L_n(\cdot) \in A}\geq
\sup_{t\in T}\pp{\frac{1}{n} L_n(t) \geq \zeta(t)}.
\end{align}
Theorem \ref{thm:crossing} indicates that the inequality in
(\ref{eq:unionBound}) is actually tight (under the conditions
stated). Informally, this means that the contribution of the
maximizing $t$ in the right-hand side of (\ref{eq:unionBound}),
say $t^\star$, dominates the contributions of the other time
epochs as $n$ grows large. This essentially says that {\it given}
that the rare event under consideration occurs, with overwhelming
probability it happens at time $t^\star.$ \hfill$\diamondsuit$
\end{rem}

As is clear from the statement of Theorem \ref{thm:crossing}, two
assumptions are needed to prove the claim; we now briefly comment
on the role played by these.

\begin{rem}
Assumption (\ref{eq:thm:crossing:unicity}) is needed to make sure
that there is not a time epoch $\bar t$, different from $t^\star$,
having a contribution of the same order as $t^\star.$ It can be
verified from our proof that if the uniqueness assumption is not
met, the probability under consideration remains asymptotically
proportional to $e^{-nI}/\sqrt{n}$, but we lack a clean expression
for the proportionality constant.

Assumption (\ref{eq:thm:crossing:limit}) has to be imposed to make
sure that the contribution of the `upper tail', that is, time
epochs $t\in \{t^\star+1, t^\star+2,\ldots\}$, can be neglected;
more formally, we should have
\begin{align}
\nt \pp{\exists t\in\{t^\star+1,t^\star+2,\ldots\}: \frac{1}{n}
L_n(t) \ge \zeta(t)} &= o \left(\pp{\frac{1}{n} L_n(\cdot) \in
A}\right).
\end{align}
In order to achieve this, the probability that the  normalized
loss process exceeds $\zeta$ for large $t$ should be sufficiently
small. \hfill$\diamondsuit$
\end{rem}

\begin{rem}
We now comment on what Assumption (\ref{eq:thm:crossing:limit})
means. Clearly,
\begin{align}
\nt \Lambda_{UZ(t)}(\theta) &= \log{\pp{\tau \le t} \ee {e^{\theta
U}}+ \pp{\tau > t}}\leq \log \ee{ e^{\theta U}},
\end{align}
as $t$ grows, $\theta\ge 0$; the limiting value is actually $\log
\ee {e^{\theta U}}$ if $\tau$ is non-defective. This entails that
\begin{align}
\nt I_{UZ}(t) = \Lambda_{UZ(t^\star)}^\star(\zeta(t)) \ge
\Lambda_U^\star(\zeta(t)) = \sup_\theta\left(\theta \zeta(t) -
\log \ee{ e^{\theta U}}\right).
\end{align}
We observe that Assumption (\ref{eq:thm:crossing:limit}) is
fulfilled if $\liminf_{t\to\infty} \Lambda_U^\star(\zeta(t))/\log
t >0$, which turns out to be valid under extremely mild
conditions. Indeed, relying on Lemma \ref{lem:app:DZ2.2.20}, we
have that in great generality it holds that
$\Lambda_U^\star(x)/x\to\infty$ as $x\to\infty$. Then clearly any
$\zeta(t)$, for which $\liminf_t\zeta(t)/\log t>0$, satisfies
Assumption (\ref{eq:thm:crossing:limit}), since
\begin{align}
\nt \liminf_{t\to\infty} \frac{\Lambda_U^\star(\zeta(t))}{\log t}
&= \liminf_{t\to\infty}
\frac{\Lambda_U^\star(\zeta(t))}{\zeta(t)}\frac{\zeta(t)}{\log t}.
\end{align}
Alternatively, if $U$ is chosen distributed exponentially with
mean $\lambda$ (which does not satisfy the conditions of Lemma
\ref{lem:app:DZ2.2.20}), then $\Lambda_U^\star(t)=\lambda t -1
-\log(\lambda t)$, such that we have that
\begin{align}
\nt \liminf_{t\to\infty} \frac{I_U(\log t )}{\log t} = \lambda >0.
\end{align}
Barrier functions $\zeta$ that grow at a rate slower than $\log
t$, such as $\log\log t$, are in this setting clearly not allowed.
\hfill$\diamondsuit$
\end{rem}

\beginproofof{Theorem \ref{thm:crossing}}

We start by rewriting the probability of interest as
\begin{align}
\nt \pp{\frac{1}{n}L_n(\cdot)\in A} &=\pp{\exists t\in T:\
\frac{L_n(t)}{n} \geq \zeta(t)}.
\end{align}
For an arbitrary point $k$ in $T$ we have
\begin{equation}
\label{eq:thm:crossing:split} \pp{\exists t\in T:\
\frac{L_n(t)}{n} \geq \zeta(t)} \leq \pp{\exists t\leq k:\
\frac{L_n(t)}{n} \geq \zeta(t)}+\pp{\exists t>k:\ \frac{L_n(t)}{n}
\geq \zeta(t)}.
\end{equation}
We first focus on the second part in
(\ref{eq:thm:crossing:split}). We can bound this by
\begin{align}
\nt
\pp{\exists t>k:\ \frac{L_n(t)}{n} \geq \zeta(t)} &\leq \sum_{i=k+1}^\infty \pp{\frac{L_n(i)}{n} \geq \zeta(i)}\\
\nt &\leq \sum_{i=k+1}^\infty
\inf_{\theta>0}\ee{\exp\left(\theta\sum_{j=1}^n
U_jZ_j(i)\right)}e^{-n\zeta(i)\theta},
\end{align}
where the second inequality is due to the Chernoff bound
\cite[Eqn.\ (2.2.12)]{book:DemboZeitouni}. The independence
between the $U_i$ and $Z_i(t)$, together with the assumption that
the $U_i$ are i.i.d.\ and the $Z_i(t)$ are i.i.d., yields
\begin{align}
\nt
\lefteqn{\sum_{i=k+1}^\infty \inf_{\theta>0}\ee{\exp\left(\theta\sum_{j=1}^n U_jZ_j(i)\right)} e^{-n\zeta(i)\theta} = \sum_{i=k+1}^\infty \inf_{\theta>0} \prod_{j=1}^n \ee{\exp\left(\theta U_jZ_j(i)\right)}e^{-n\zeta(i)\theta}}\\
\nt &= \sum_{i=k+1}^\infty
\exx{-n\sup_{\theta>0}\left(\zeta(i)\theta -
\Lambda_{UZ(i)}(\theta\right)} =  \sum_{i=k+1}^\infty
\exx{-nI_{UZ}(\zeta(i))}.
\end{align}
By (\ref{eq:thm:crossing:limit}) we have that
\begin{align}
\nt \liminf_{t\to\infty} \frac{I_{UZ}(t)}{\log t}&=\beta,
\end{align}
for some $\beta > 0$ (possibly $\infty$). Hence there exists an
$m$ such that for all $i>m$
\begin{equation}
\label{eq:thm:crossing:ineqs} I_{UZ}(i) > \alpha\log i >
I_{UZ}(t^\star),
\end{equation}
where $\alpha = {\beta}/{2}$ (in case $\beta=\infty$, any
$0<\alpha<\infty$ suffices) and $t^\star$ defined in
(\ref{eq:thm:crossing:unicity}). Choosing $k = m$, we obtain by
using the first inequality in (\ref{eq:thm:crossing:ineqs}) for
$n> 1/\alpha$
\begin{align}
\nt \sum_{i=m+1}^\infty \exx{-nI_{UZ}(\zeta(i))} &\leq
\sum_{i=m+1}^\infty \exx{-n\alpha\log i } \leq
\frac{1}{n\alpha-1}\exx{(-n\alpha+1)\log m },
\end{align}
where the last inequality trivially follows by bounding the
summation (from above) by an appropriate integral. Next we
multiply and divide this by $\pp{L_n(t^\star)/n>\zeta(t^\star)}$
and we apply the Bahadur-Rao theorem, which results in
\begin{align}
\nt\lefteqn{ \frac{1}{n\alpha-1}e^{(-n\alpha+1)\log m} =
\frac{1}{n\alpha-1}e^{(-n\alpha+1)\log m}\
\frac{\pp{L_n(t^\star)/n>\zeta(t^\star)}}{\pp{ L_n(t^\star)/n>\zeta(t^\star)}}}\\
\nt&= \pp{\frac{1}{n}L_n(t^\star)>\zeta(t^\star) } \
\frac{m\sqrt{n}\,C^\star}{n\alpha-1}\
\left(1+O\left(\frac{1}{n}\right)\right) \ e^{-n\left(\alpha\log m
-I_{UZ}(t^\star) \right)}.
\end{align}
The second inequality in (\ref{eq:thm:crossing:ineqs}) yields
$\alpha\log m -I_{UZ}(t^\star) > \delta$, for some $\delta > 0$.
Applying this inequality, we see that this bounds the second term
in (\ref{eq:thm:crossing:split}), in the sense that as
$n\to\infty$,
\begin{align} \nt \left. \pp{\exists t>k:\ \frac{L_n(t)}{n} \geq \zeta(t)}\right/\pp{\frac{1}{n} L_n(t^\star)>\zeta(t^\star)}\to 0.
\end{align}
To complete the proof we need to bound the first term of
(\ref{eq:thm:crossing:split}), where we use that $k=m$. For this
we again use the Bahadur-Rao theorem. Next to this theorem we use
the uniqueness of $t^\star$, which implies that for $i\leq m$ and
$i\neq t^\star$ there exists an $\epsilon^\star>0$, such that
\begin{align}
\nt I_{UZ}(t^\star)+\epsilon^\star\leq I_{UZ}(i).
\end{align}
This observation yields, with $\sigma_i$ such that
$\Lambda_{UZ(i)}'(\sigma_i)=\zeta(i)$,
\begin{align}
\nt
\lefteqn{\pp{\exists t\leq m:\ \frac{L_n(t)}{n} \geq \zeta(t)} \leq \sum_{i=1}^m\pp{  \frac{L_n(i)}{n} \geq \zeta(i)}} \\
\nt &\leq\pp{  \frac{1}{n}L_n(t^\star)>\zeta(t^\star)}\ \left( 1+O \left( \frac{1}{n} \right)\right)\ \left(\sum_{i=1}^m \frac{C^\star}{C_{UZ(i),\zeta(i)}}\ \frac{e^{  -n I_{UZ}(t_i)}}{e^{  -n I_{UZ}(t^\star) }}\right) \\
\nt &\leq\pp{ \frac{1}{n}L_n(t^\star)>\zeta(t^\star)}\ \left( 1+O \left( \frac{1}{n} \right)\right) \ \left(  1+ m\times\max_{i=1,\ldots,m}  \left(   \frac{C^\star}{C_{UZ(i),\zeta(i)}}   \right)  e^{ - n \varepsilon^\star}\right)\\
\nt &=\pp{ \frac{1}{n}L_n(t^\star)>\zeta(t^\star)}\ \left(  1+O
\left( \frac{1}{n} \right)\right)\ \left(  1+ O  \left(    e^{ - n
\varepsilon^\star}  \right)\right)
\end{align}
Combining the above findings, we  observe
\begin{align}
\nt \pp{  \exists t\in T:\ \frac{L_n(t)}{n}\ge \zeta(t)} \leq
\pp{\frac{L_n(t^\star)}{n}\ge\zeta(t^\star)}\left(  1+O \left(
\frac{1}{n} \right)\right).
\end{align}
Together with the trivial bound
\begin{align}
\nt \pp{  \exists t\in T:\ \frac{L_n(t)}{n}\ge \zeta(t)} \geq
\pp{\frac{L_n(t^\star)}{n} \ge \zeta(t^\star)},
\end{align}
this yields
\begin{align}
\nt \pp{\exists t\in T:\ \frac{L_n(t)}{n}\ge \zeta(t)} =
\pp{\frac{L_n(t^\star)}{n}>\zeta(t^\star)}\left(  1+O \left(
\frac{1}{n} \right)\right).
\end{align}
Applying the Bahadur-Rao theorem to the right-hand side of the
previous display yields the desired result.
\end{proof}

\subsection{Large Increments of the Loss Process}
In the previous section we identified the asymptotic behavior of
the probability that at some point in time the normalized loss
process $L_n(\cdot)/n$ exceeds a certain level. We can carry out a
similar procedure to obtain insight in the large deviations of the
{\it increments} of the loss process. Here we consider times where
the increment of the loss between time $s$ and $t$ exceeds a
threshold $\xi(s,t)$. More precisely, we consider the event
\begin{align}
\label{eq:crossingIntervalEvent} A:= \left\{f:T\times
T\to\R_0^+\vert\ \exists s<t\in T:\ f(s,t) \geq \xi(s,t)\right\}.
\end{align}
Being able to deal with events of this type, we can for instance
analyze the likelihood of the occurrence of a large loss during a
short period; we remark that with the event (\ref{eq:crossingSet})
from the previous subsection, one cannot distinguish the cases
where the loss is zero for all times before $t$ and $x>\zeta(t)$
at time $t$, and the case where the loss is just below the level
$\zeta$ for all times before time $t$ and then ends up at $x$ at
time $t$. Clearly, events of the (\ref{eq:crossingIntervalEvent})
make it possible to distinguish between such paths.

In order to avoid trivial results, we impose a condition similar
to (\ref{eq:interestingLevel}), namely \begin{align}
\label{eq:interestingLevelIncrement} \xi(s,t) &> \ee{U}\
\left(F_t-F_s\right),
\end{align}
for all $s<t$. The law of large numbers entails that for functions
$\xi$ that do not satisfy this condition, the probability under
consideration does not correspond to a rare event.

A similar probability has been considered in
\cite{article:DemboDeuschelDuffie2004},  where the authors derive
the logarithmic asymptotic behavior of the probability that the
increment of the loss, for some $s<t$, in a bounded interval
exceeds a thresholds that depends only on $t-s$. In contrast, our
approach uses a more flexible threshold, which depends on both
times $s$ and $t$, and in addition we derive the {\it exact}
asymptotic behavior of this probability.

\begin{thm}
\label{thm:crossingInterval} Assume that
\begin{align}
\label{eq:thm:crossingInterval:uniq}\mbox{there is a unique
$s^\star<t^\star\in T$ such that}\:\:  I_{UZ}(s^\star,t^\star)&=
\min_{s<t}I_{UZ}(s,t),
\end{align}
and that
\begin{align}
\label{eq:thm:crossingInterval:limit} \inf_{s\in
T}\liminf_{t\to\infty}\frac{I_{UZ}(s,t)}{\log t}>0,
\end{align}
where $I_{UZ}(s,t) = \sup_\theta \left(\theta\xi(s,t) -
\Lambda_{U(Z(t)-Z(s))} (\theta)\right) =
\Lambda_{U(Z(t)-Z(s))}^\star (\xi(s,t)).$ Then
\begin{align}
\label{eq:thm:crossingInterval:result}
\pp{\frac{1}{n} L_n(\cdot)\in A} =
\frac{e^{-nI_{UZ}(s^\star,t^\star)}C^\star}{\sqrt{n}}\order,
\end{align}
for $A$ as in (\ref{eq:interestingLevelIncrement}) and
$\sigma^\star$ is such that
$\Lambda_{U(Z(t^\star)-Z(s^\star))}'(\sigma^\star) =
\xi(s^\star,t^\star)$. The constant $C^\star$ follows from the
Bahadur-Rao theorem (Theorem \ref{thm:bahadurRao}), with
$C^\star=C_{U(Z(t^\star)-Z(s^\star)),\ \xi(s^\star,t^\star)}$.
\end{thm}

\begin{rem}
A first glance at Theorem \ref{thm:crossingInterval} tells us the
obtained result is very similar to the result of Theorem
\ref{thm:crossing}. The second condition, i.e., Inequality
(\ref{eq:thm:crossingInterval:limit}), however, seems to be more
restrictive than the corresponding condition, i.e., Inequality
(\ref{eq:thm:crossing:limit}), due to the infimum over $s$. This
assumption has to make sure that the `upper tail' is negligible
for any $s$. In the previous subsection we have seen that, under
mild restrictions, the upper tail can be safely ignored when the
barrier function grows at a rate of at least $\log t$. We can
extend this claim to our new setting of large increments, as
follows.

First note that
\begin{align}
\nt \inf_{s\in T}\liminf_{t\to\infty}\frac{I_{UZ}(s,t)}{\log
t}&\geq \inf_{s\in
T}\liminf_{t\to\infty}\frac{\Lambda_U^\star(\xi(s,t))}{\log t}.
\end{align}

Then consider thresholds that, next to condition
(\ref{eq:interestingLevelIncrement}), satisfy that for all $s$
\begin{align}
\label{eq:sufficientThresholdInterval}
\liminf_{t\to\infty}\frac{\xi(s,t)}{\log t} &>0.
\end{align}
Then, under the conditions of Lemma \ref{lem:app:DZ2.2.20}, we
have that
\begin{align}
\label{eq:toInfForAllS}
\liminf_{t\to\infty}\frac{\Lambda_U^\star(\xi(s,t))}{\log t} &=
\liminf_{t\to\infty}\frac{\Lambda_U^\star(\xi(s,t))}{\xi(s,t)}\frac{\xi(s,t)}{\log
t}=\infty,
\end{align}
since the second term remains positive by
(\ref{eq:sufficientThresholdInterval}) and the first term tends to
infinity by 
Lemma \ref{lem:app:DZ2.2.20}. Having established
(\ref{eq:toInfForAllS}) for all $s$, it is clear that
(\ref{eq:thm:crossingInterval:limit}) is satisfied.

The sufficient condition (\ref{eq:sufficientThresholdInterval})
shows that the range of admissible barrier functions is quite
substantial, and, importantly, imposing
(\ref{eq:thm:crossingInterval:limit}) is not as restrictive as it
seems at first glance. \hfill$\diamondsuit$
\end{rem}

\beginproofof{Theorem \ref{thm:crossingInterval}}

The proof of this theorem is very similar to that of Theorem
\ref{thm:crossing}. Therefore we only sketch the proof here.

As before, the probability of interest is split up into a `front
part' and `tail part'. The tail part can be bounded using
Assumption (\ref{eq:thm:crossingInterval:limit}); this is done
analogously to the way Assumption (\ref{eq:thm:crossing:limit})
was used in the proof of Theorem \ref{thm:crossing}. The
uniqueness assumption (\ref{eq:thm:crossingInterval:uniq}) then
shows that the probability of interest is asymptotically equal to
the probability that the increment between time $s^\star$ and
$t^\star$ exceeds $\xi(s^\star,t^\star)$; this is an application
of the Bahadur-Rao theorem. Another application of the Bahadur-Rao
theorem to the probability that the increment between time
$s^\star$ and $t^\star$ exceeds $\xi(s^\star,t^\star)$ yields the
result.
\end{proof}

\section{Discussion and Concluding Remarks}
\label{sec:discussion} In this paper we have established a number
of results with respect to the asymptotic behavior of the
distribution of the loss process. In this section we discuss some
of the assumptions in more detail and we consider extensions of
the results that we have derived.

\subsection{Extensions of the Sample-Path LDP}
The first part of our work, Section \ref{sec:samplePathLDP}, was
devoted to establishing a sample-path large deviation principle on
a finite time grid. Here we modelled the loss process as the sum
of i.i.d.\ loss amounts multiplied by i.i.d.\ default indicators.
>From a practical point of view one can argue that the assumptions
underlying our model are not always realistic. In particular, the
random properties of the obligors cannot always be assumed
independent. In addition, the assumption that all obligors behave
in an i.i.d.\ fashion will not necessarily hold in practice. Both
shortcomings can be dealt with, however, by adapting the model
slightly.

A common way to introduce dependence, taken from
\cite{article:DemboDeuschelDuffie2004}, is by supposing that there
is a `macro-environmental' variable $Y$ to which all obligors
react, but conditional on which the loss epochs and loss amounts
are independent. First observe that our results are then valid for
any specific realization $y$ of $Y$. Denoting the exponential
decay rate by $r_y$, i.e.,
\[\lim_{n\to\infty}\frac{1}{n}\log \pp{\left.\frac{1}{n}L_n(\cdot)\in A\:\right\vert\:Y=y}=r_y,\]
the {\it unconditional} decay rate is just the maximum over the
$r_y$; this is trivial to prove if $Y$ can attain values in a
finite set only. A detailed treatment of this is beyond the scope
of this paper.

The assumption that all obligors have the same distribution can be
relaxed to the case where we assume that  there are $m$ different
classes of obligors (for instance determined by their
defaultratings). We further assume that each class $i$ makes up a
fraction $a_i$ of the entire portfolio. Then we can extend the LDP
of Theorem \ref{subsec:samplePathLDPmainResult} to a more general
one, by splitting up the loss process into $m$ loss processes,
each corresponding to a class. Conditioning on the realizations of
these processes, we can derive the following rate function:
\begin{align}
\label{eq:multiTypeLDP} I_{U,p,m}(x)&:= \inf_{\phi\in\Phi^m}
\inf_{v\in V_x} \sum_{j=1}^m \sum_{i=1}^N a_i \phi_i^j
\left(\lgg{\frac{\phi_i^j}{p_i^j}}+\Lambda_U^\star\left(\frac{
v_i^j}{a_i \phi_i^j}\right)\right),
\end{align}
where $V_x=\{\left.v\in \R_+^{m\times N}\right\vert\ \sum_{j=1}^m
v_i^j = \Delta x_i \mbox{ for all } i\leq N\}$, and $\Phi^m$ is
the Cartesian product $\Phi\times\ldots\times\Phi$ ($m$ times),
with $\Phi$ as in (\ref{eq:defPhiN}). The optimization over the
set $V_x$ follows directly from conditioning on the realizations
of the per-class loss processes. We leave out the formal
derivation of this result; this multi-class case is notationally
considerably more involved than the single-class case, but
essentially all steps carry over.

In our sample-path LDP we assumed that defaults can only occur on
a finite grid. While this assumption is justifiable from a
practical point of view, an interesting mathematical question is
whether it can be relaxed. In self-evident notation, one would
expect that the rate function
\begin{align}
\nt  I_{U,p,\infty}(x)&:=
\inf_{\phi\in\Phi_\infty}\sum_{i=1}^\infty
\phi_i\left(\lgg{\frac{\phi_i}{p_i}}
+\Lambda_U^\star\left(\frac{\Delta x_i}{\phi_i}\right)\right).
\end{align}
It can be checked, however, that the argumentation used in the
proof of Theorem \ref{thm:resultFiniteGrid} does not work; in
particular, the choice of a suitable topology plays an important
role.

If losses can occur on a continuous entire interval, i.e.,
$[0,N]$, we expect, for a nondecreasing and differentiable path
$x$, the rate function
\begin{align}
\label{eq:con:intervalLDP} I_{U,p,[0,N]}(x)&:=
\inf_{\phi\in\mathscr{M}} \int_0^N
\phi(t)\left(\lgg{\frac{\phi(t)}{p(t)}}+
\Lambda_U^\star\left(\frac{x'(t)}{\phi(t)}\right)\right)\dd t,
\end{align}
where $\mathscr{M}$ is the space of all densities on $[0,N]$  and
$p$ the density of the default time $\tau$.  One can easily guess
the validity of (\ref{eq:con:intervalLDP}) from
(\ref{eq:defRateFctFiniteGrid}) by using Riemann sums to
approximate the integral. A formal proof, however, requires
techniques that are essentially different from  the ones used to
establish Theorem \ref{thm:resultFiniteGrid}, and therefore we
leave this for future research.

\subsection{Extensions of the Exact Asymptotics}

In the second part of the paper, i.e., Section
\ref{sec:exactResults}, we have derived the exact asymptotic
behavior for two special events. First we showed that, under
certain conditions, the probability that the loss process exceeds
a certain time-dependent level, is asymptotically equal to the
probability that the process exceeds this level at the `most
likely' time $t^\star$. The exact asymptotics of this probability
are obtained by applying the Bahadur-Rao theorem. A similar result
has been obtained  for an event related to the {\it increment} of
the loss process. One could think of refining the logarithmic
asymptotics, as developed in Section \ref{sec:samplePathLDP}, to
exact asymptotics. Note, however, that this is far from
straightforward, as for general sets these asymptotics do not
necessarily coincide with those of a univariate random variable,
cf.\ \cite{MMNU}.

\bibliographystyle{amsalpha}

\appendix

\section{Background Results}
In this section we state a number of definitions and results,
taken from \cite{book:DemboZeitouni}, which are used in the proofs
in this paper.

\begin{thm}[Cram\'er]
\label{thm:Cramer} Let $X_i$ be i.i.d.\ real valued random
variables with all exponential moments finite and let $\mu_n$ be
the law of the average $S_n= \sum_{i=1}^nX_i/n$. Then the sequence
$\{\mu_n\}$ satisfies an LDP with rate function
$\Lambda^\star(\cdot)$, where $\Lambda^\star$ is the
Fenchel-Legendre transform of the $X_i$.
\end{thm}
\beginproof
See for example \cite[Thm.\ 2.2.3]{book:DemboZeitouni}.
\end{proof}

\begin{df}
\label{df:expEquiv} We say that two families of measures
$\{\mu_n\}$ and $\{\nu_n\}$ on a metric space $(\X,d)$ are
exponentially equivalent if there exist two families of
$\X$-valued random variables $\{Y_n\}$ and $\{Z_n\}$ with marginal
distributions $\{\mu_n\}$ and $\{\nu_n\}$, respectively, such that
for all $\delta>0$
\begin{align}
\nt
\limsup_{n\to\infty}\frac{1}{n}\log{\pp{d\left(X_n,Y_n\right)\geq\delta}}=-\infty.
\end{align}

\end{df}

\begin{lemma}
\label{lem:app:DZ1.2.15} For  every triangular array $a_n^i\geq
0$, $n\geq 1$, $1\leq i\leq n$,
\begin{align}
\nt \limsup_{n\to\infty}\frac{1}{n}\log{\sum_{i=1}^n a_n^i} &=
\limsup_{n\to\infty}\max_{i=1,\ldots,n}\frac{1}{n}\log{a_n^i}.
\end{align}
\end{lemma}
\beginproof
Elementary, but also a direct consequence of \cite[Lemma
1.2.15]{book:DemboZeitouni}.
\end{proof}

\begin{lemma}
\label{lem:app:DZ2.2.20} Let $\Lambda(\theta)<\infty$ for all
$\theta\in\R$, then
\begin{align}
\nt \lim_{\vert x \vert \to \infty}\frac{\Lambda^\star (x)}{\vert
x\vert } = \infty.
\end{align}
\end{lemma}
\beginproof
This result is a part of \cite[Lemma 2.2.20]{book:DemboZeitouni}.
\end{proof}

\begin{lemma}
\label{lem:app:DZ2.1.9} Let $K_{n,i}$ be defined as
$K_{n,j}:=\#\{i\in\{1,\ldots,n\}\:\vert\: \tau_i=j\}.$ Then for
any vector $k\in\N^N$, such that $\sum_{i=1}^N k_i=n$, we have
that
\begin{align}
\nt (n+1)^{-N}\exx{-nH(k\:\vert\: p)}\leq  \pp{K_n = k} \leq
\exx{-nH(k\:\vert\: p)},
\end{align}
where
\begin{align}
\nt H(k\:\vert\: p) &=\sum_{i=1}^N \frac{k_i}{n}\lgg{\frac{k_i}{n
p_i}},
\end{align}
and $p_i$ as defined in $(\ref{eq:defLossDist}).$
\end{lemma}
\beginproof
See \cite[Lemma 2.1.9]{book:DemboZeitouni}.
\end{proof}

\begin{lemma}
\label{lem:app:DZ5.1.8} Define
\begin{align}
\nt Z_n(t):=\frac{1}{n}\sum_{i=1}^{[nt]}X_i,\ 0\leq t\leq 1,
\end{align}
for an i.i.d.\ sequence of $\R^d$  valued random variables $X_i$.
Let $\mu_n$ denote the law of $Z_n(\cdot)$ in $L_\infty([0,1])$.
For any discretization $J=\{0<t_1<\ldots<t_{\vert J\vert}\leq 1\}$
and any $f:[0,1]\to\R^d$, let $p_J(f)$ denote the vector
$(f(t_i))_{i=1}^{\vert J\vert}\in (\R^d)^{\vert J \vert}$. Then
the sequence of laws $\{\mu_n\circ p_J^{-1}\}$ satisfies the LDP
in $(\R^d)^{\vert j \vert}$ with the good rate function
\begin{align} \nt I_J(z)=\sum_{i=1}^{\vert J
\vert}(t_i-t_{i-1})\Lambda^\star\left(\frac{z_i-z_{i-1}}{t_i-t_{i-1}}\right),
\end{align}
where $\Lambda^\star$ is the Fenchel-Legendre transform of $X_1$.
\end{lemma}
\beginproof
See \cite[Lemma 5.1.8]{book:DemboZeitouni}. This lemma is one of
the key steps in proving  $\MOG$'s theorem, which provides a
sample-path LDP for $Z_n(\cdot)$ on a bounded interval.
\end{proof}

\begin{thm}
\label{thm:app:DZ4.2.13} If an LDP with a good rate function
$I(\cdot)$ holds for the probability measures $\{\mu_n\}$, which
are exponentially equivalent to $\{\nu_n\}$, then the same LDP
holds for $\{\nu_n\}$.
\end{thm}
\beginproof
See  \cite[Thm.\ 4.2.13]{book:DemboZeitouni}.
\end{proof}

\begin{thm}[Bahadur-Rao]
\label{thm:bahadurRao} Let $X_i$ be an i.i.d.\ real-valued
sequence random variables. Then we have
\begin{align}
\nt \pp{\frac{1}{n}\sum_{i=1}^n X_i \ge q} = \frac{
e^{-n\Lambda_X^\star(q)}C_{X,q}}{\sqrt{n}}\order.
\end{align}
The constant $C_{X,q}$ depends on the type of distribution of
$X_1$, as specified by the following two cases.
\begin{itemize}
\item[(i)] The law of $X_1$ is lattice, i.e. for some $x_0$, $d$, the random variable $(X_1-x_0)/d$ is  (a.s.) an integer number, and $d$ is the largest number with this property. Under the additional condition $0<\pp{X_1=q}<1$, the constant $C_{X,q}$ is given by
\[C_{X,q}=\frac{d}{(1-e^{-\sigma d})\sqrt{2\pi \Lambda_X''(\sigma)}},\]
where $\sigma$ satisfies $\Lambda_X'(\sigma)=q$.
\item[(ii)] If the law of $X_1$ is non-lattice, the constant $C_{X,q}$ is given by
\[C_{X,q}=\frac{1}{\sigma\sqrt{2\pi \Lambda_X''(\sigma)}},\]
with $\sigma$  as in case \emph{(i)}.
\end{itemize}

\end{thm}
\beginproof
We refer to \cite{article:BahadurRao1960} or \cite[Thm.\
3.7.4]{book:DemboZeitouni} for the proof of this result.
\end{proof}

\vspace{1cm}

\subsection*{Aknowledgment} {\small VL would like to thank ABN A{\sc mro} bank for providing financial support. Part of this work was carried out while MM was at Stanford University, US. The authors are indebted to E.J.\ Balder (Utrecht University, the Netherlands) for pointing out to us the relevance of epi-convergence to our research.}
\end{document}